\documentclass{article}
\usepackage[T2A]{fontenc}
\usepackage[utf8]{inputenc}
\usepackage[tbtags]{amsmath}
\usepackage{amsfonts,amssymb,mathrsfs,amscd,comment}

\begin{document}


\title{\bf\large\MakeUppercase{%
On coset $n$-valued topological groups \\ on $S^3$ and $\mathbb{R}P^3$
}}

\author{Dmitry\,V.\,Gugnin}
\date{}
\maketitle
\begin{flushright}
\large {\it To my son Vanya}
\end{flushright}

\makeatletter
\renewcommand{\@makefnmark}{}
\makeatother
\footnotetext{This work was supported by the Russian Science Foundation under grant no.~20-11-19998, https://rscf.ru/en/project/20-11-19998/ . }    

Recall the notion of an $n$-valued topological group (see an overview \cite{Buch2006}). 

\smallskip

\textbf{Definition 1.} {\it A Hausdorff path-connected topological space $X$ with a base point $e\in X$, together with an $n$-valued continuous multiplication $\mu\colon X\times X \to \mathrm{Sym}^nX = X^n/S_n$ and a continuous involution $\mathrm{inv}\colon X \to X, \mathrm{inv}(e) = e$, is called an $n$-valued topological group, if the following conditions hold true:\\
(1) $\mu(x,\mu(y,z)) = \mu(\mu(x,y),z)\in \mathrm{Sym}^{n^2}X$ for all $x,y,z\in X$;\\
(2) $\mu(x,e) = \mu(e,x) = [x,x,\ldots,x]$ for all $x\in X$;\\
(3) $\mu(x, \mathrm{inv}(x))\ni e$ и $\mu(\mathrm{inv}(x), x)\ni e$ for all $x\in X$.}

\smallskip

We will use the following construction of a {\it coset} $n$-valued groups (see \cite{Buch2006}).
\smallskip

Let us take an arbitrary compact connected Lie group $W$ and any subgroup $G\subset \mathrm{Aut}(W)$ of order $n$ in its automorphism group. Denote by $X$ the quotient space $W/G$ and by $\pi\colon W \to X$ the canonic projection. Then the space $X$ can be endowed with a natural structure of an $n$-valued topological group with the identity $e=\pi(e_W)$, the inverse map $\mathrm{inv} \colon X \to X, \mathrm{inv}(\pi(w)) = \pi(w^{-1})$ and multiplication  $\mu(\pi(a),\pi(b)) = [\pi(a g_1(b)), \pi(a g_2(b)), \ldots, \pi(a g_n(b))]$ for all $a,b\in W$, where $G=\{g_1,g_2,\ldots, g_n\}$. 

\smallskip

In dimension 3 there are only three compact connected Lie groups: $T^3, Sp(1)$ и $SO(3)$. The aim of the research is to describe all coset $n$-valued topological groups, arising from $Sp(1)$ and $SO(3)$.

\smallskip

It is known that any automorphism $\varphi\colon SO(3) \to SO(3)$ is an inner conjugation by an element $g\in SO(3)$, which is 1-1 correspondence to $\varphi$. Therefore the group $\mathrm{Aut}(SO(3))$ coincides with $SO(3)$.   Also the classical isomorphism  $Sp(1)/\{1,-1\} \cong SO(3)$ implies that any automorphism $\psi\colon Sp(1) \to Sp(1)$ is an inner conjugation by an element $\pm q\in Sp(1)$. Hence the group $\mathrm{Aut}(Sp(1))$ coincides again with $SO(3)$.

\smallskip

Let us use the well known classification of finite subgroups in SO(3).

\smallskip

\noindent (1) $C_{n}$, a cyclic group, generated by a $n$-fold rotation about a line;  \\
(2) $D_{m}$, $n=2m$, a dihedral group, generated by an $n$-fold rotation about a line, and a reflection in a line
(half-turn) which is orthogonal to the first line; \\
(3) $T$, $n=12$, the group of orientation-preserving symmetries of a regular tetrahedron; \\
(4) $O$, $n=24$, the group of orientation-preserving symmetries of a cube; \\
(5) $I$, $n=60$, the group of orientation-preserving symmetries of a regular icosahedron.

\smallskip
The following fact is known: 

\smallskip

\textbf{Fact 1.} { \it For an arbitrary smooth connected orientable manifold $M^m, m=2,3$ and any finite group $G$, acting on $M^m$ smoothly, effectively and orientation-preserving, the orbit space $X=M^m/G$ is a topological orientable $m$-dimensional  manifold.}

\smallskip

Therefore, if a coset topological group  $X$ is derived from $Sp(1)$ or $SO(3)$, then $X$ is an orientable connected compact topological 3-manifold.
\smallskip

\textbf{Theorem 1.}  {\it Let $W=Sp(1)$ and $G\subset \mathrm{Aut}(Sp(1)) = SO(3)$ --- an arbitrary finite subgroup of order  $n$. Then the coset $n$-valued topological group $Sp(1)/G$ is homeomorphic to $S^3$.} 

\smallskip

The case of a 2-valued coset group structure on $S^3$ was introduced by V.\,M.~Buchstaber in 1993.

\smallskip

\textbf{Theorem 2.}  {\it Let $W=SO(3)$ and $G\subset  SO(3)$ --- an arbitrary finite subgroup of order $n$. If $n$ is even, then the coset group $SO(3)/G$ is homeomorphic to $S^3$. And if $n$ is odd, then the coset group $SO(3)/G$ is homeomorphic to $\mathbb{R}P^3$.}

\smallskip

\textbf{Proof of theorem 1.} It is obvious that the action of the group $G$ preserves the real part of quaternions $x\in Sp(1)$, fix points $\pm1$, and also preserves the standard metric of the sphere $S^3=Sp(1)$. It follows that the quotient space $Sp(1)/G$ is the (unreduced) suspension over the quotient space $S^2/G$, where $S^2$  is a sphere of purely imaginary quaternions of unit length. By the fact that the action of the group $G$ is smooth and orientation-preserving and due to fact 1, the quotient space $S^2/G$ is a compact orientable surface $M^2$. As the map $S^2 \to M^2 = S^2/G$ has a nonzero degree, then the genus of a surface $M^2$ could not increase, so $M^2 = S^2$. \ \ \  $\Box$ 

\smallskip

\textbf{Proof of theorem 2.} Let us take a finite subgroup $G\subset SO(3)$ of order $n$. Denote by $2G$ its preimage under the canonic epimorphism $Sp(1) \to SO(3)$. Set $G=\{g_1,g_2,\ldots, g_n \}$ and $2G=\{\pm q_1, \pm q_2, \ldots, \pm q_n\}$. The quoteint space $SO(3)/G$ may be derived from the universal cover $\tilde{SO}(3) = Sp(1) = S^3$ by the following action of a finite group $\tilde{G} := G\times C_2, C_2 =\{1,-1\}$. Namely, $(g_i, \varepsilon) (x) := \varepsilon (q_ixq_i^{-1})$, for all $x\in Sp(1)$ and $\varepsilon \in \{1,-1\}$. The required orbit space $Sp(1)/\tilde{G}$ could be reached into two steps. Firstly, we need to get the quotient space $Sp(1)/G\cong S^3$. Secondly, we need to factorize the obtained sphere $S^3$ by an involution $\tau$, which arise from antipodal involution $x\mapsto -x$ on $Sp(1)$. It is clear, that the involution $\tau$ preserves the orientation. 

\smallskip

Classical S.Illman's theorem \cite{Illman78} states that for any smooth manifold $M^m$ and any finite group $F$, which acts smoothly on $M^m$, there exists a triangulation of $M^m$ for which the action of $F$ is simplicial. It follows that the involution $\tau\colon S^3 \to S^3$ is simplicial. 

\smallskip

F.Waldhausen' theorem \cite{Waldhausen69} states that any preserving orientation simplicial involution of the sphere $S^3$ is conjugated in the group of homeomorphisms to one of standard involutions:  $(y_1,y_2,y_3,y_4) \mapsto (-y_1,-y_2,-y_3,-y_4)$ (no fixed points), or $(y_1,y_2,y_3,y_4) \mapsto (-y_1,-y_2,y_3,y_4)$ (there are fixed points). In the first case the orbit space $S^3/\tau$ is $\mathbb{R}P^3$. In the second case $S^3/\tau \cong S^3$. 

\smallskip

So, to prove our theorem 2  it is sufficient to obtain a criteria for the involution $\tau$ to have fixed points. The desired existence of fixed points for $\tau$ is equivalent to the fact that for some $g_i, 1\le i \le n,$ the element $(g_i,-1)\in \tilde{G}$ possesses at least one fixed point $x\in Sp(1)$, i.e. a point $x\in Sp(1)$ with the condition $-q_ixq_i^{-1} = x$. 

\smallskip

Consider the equation $qxq^{-1} = -x$ on the group $Sp(1)$ for a given $q\in Sp(1)$. Since $\mathrm{Re} (qxq^{-1}) = \mathrm{Re}(x)$, then a solution $x$ of our equation could be only a purely imaginary quaternion of unit length. It is known, that the epimorphism $Sp(1) \to SO(3)$ maps a quaternion $q\in Sp(1)$ onto the element $g\in SO(3)$ such that $g(x) = qxq^{-1}$, where $x$ is an arbitrary purely imaginary quaternion of unit length. As any element $g\in SO(3)$ is a rotation about a line by some angle, then the equation $g(x) = -x$ has a solution iff $g$ is a rotation by  the angle $\pi$, i.e. $g^2 = e$. Since $g\in G$, then the sought condition is equivalent to the order of the group $G$ be even. \ \ \  $\Box$ 

\smallskip

The author is deeply grateful to the corresponding member of RAS, Professor V.\,M.~Buchstaber for fruitful discussions and an instant interest to the research, and to Associate Professor I.\,Yu.~Limonchenko for reading the manuscript and making valuable comments.

{\bf D.\,V.~Gugnin}

Steklov Mathematical Institute of Russian Academy of Sciences, Moscow, Russia 

{\it E-mail}: dmitry-gugnin@yandex.ru

\end{document}